\documentclass[12pt]{article}
\usepackage{float}
\usepackage[english]{babel}
\usepackage{amssymb, latexsym, amsmath, amsthm}
\usepackage[left=2cm, right=2cm, top=2cm]{geometry}
\usepackage{float}
\usepackage{tikz-cd}
\usepackage{caption}
\usepackage{float}
\usepackage{subcaption}
\usepackage{comment}
\usepackage{blindtext}
\usepackage{titlesec}
\titleformat{\section}
   {\normalfont\fontsize{12pt}{14pt}\selectfont\bfseries}{\thesection}{1em}{\MakeUppercase}
\usepackage{sectsty}
\sectionfont{\centering}
\usepackage{bbm}
\usepackage{longtable}
\usepackage{fullpage} 
\usepackage{parskip} 
\usepackage{hyperref}

\usepackage{enumerate}
\usepackage{mathtools}
\usepackage{xcolor}
\hypersetup{backref=true,       
    pagebackref=true,               
    hyperindex=true,                
    colorlinks=true,                
    breaklinks=true,                
    urlcolor= black,                
    linkcolor= black,                
    bookmarks=true,                 
    bookmarksopen=false,
    filecolor=black,
    citecolor=black,
    linkbordercolor=red
}

\newcommand{\F}{\mathbbm{F}}
\newcommand{\compconj}[1]{%
  \overline{#1}%
}
\theoremstyle{plain}
\newtheorem{thm}{Theorem}[section]
\newtheorem{prop}[thm]{Proposition}
\newtheorem{lemma}[thm]{Lemma}

\theoremstyle{definition}
\newtheorem{defn}[thm]{Definition}

\newtheorem{remark}[thm]{Remark}

\usepackage{url}
\usepackage{geometry}
\usepackage{ragged2e}
\AtBeginDocument{\hypersetup{pdfborder={0 0 1}}}
\begin{document}
\title{Zagier's weight $3/2$ mock modular form}
\author{Ajit Bhand and Ranveer Kumar Singh}
\date{}
\maketitle
\begin{abstract}
\noindent  Mock modular forms have their origins in Ramanujan's pioneering work on mock theta functions. 
In a 1975 paper, Zagier proved certain transformation properties of the generating function of the Hurwitz class numbers $H(n)$ for the discriminant $(-n)$. In the modern framework, these results show that the generating function of $H(n)$ is a mock modular form of weight 3/2 with the theta function being the shadow. In this expository paper, we provide a detailed proof of Zagier's result.
\end{abstract}
\section{Introduction}
Modular forms and their various generalizations are ubiquitous in mathematics. They enjoy diverse and often surprising applications in many areas of mathematics, including number theory, representation theory, discrete geometry \cite{wiles1, wiles2, borcherds, MV,HC} as well as theoretical physics \cite{DMZ,LS,AS,SA}. Mock modular forms have their origin in the last deathbed letter of S. Ramanujan to  G.H. Hardy written in 1920 in which he listed 17 functions with ``modular like" properties, and which he called \emph{mock theta functions} \cite{SR}. The transformation properties of Ramanujan's mock theta functions were not understood completely until 2002, when Zwegers, in his PhD. thesis, came up with a systematic framework to study mock theta functions \cite{zw}. Zwegers showed that all of Ramanujan's mock theta functions are mock modular forms of weight $1/2.$ Over the last several years, mock modular forms have become a central topic in number theory with applications ranging from geometry, topology to black hole physics (see \cite{F} for a nice review).


Modular forms are intimately connected to another central object in number theory - class numbers of quadratic fields. Roughly speaking, class numbers count the number of $\mathrm{SL}(2,\mathbb{Z})$-equivalence classes of binary quadratic forms with fixed discriminant or equivalently the index of principal fractional ideals in the group of fractional ideals of an imaginary quadratic field with given discriminant (see \cite{trifkovic} for details). The generating functions of certain generalisations of class numbers are known to be modular forms of half-integral weight \cite{Co} (see section \ref{sec 2}). In a special case these generalised class numbers reduce to the familiar Hurwitz class numbers (see Remark \ref{remark 2.5}) which count the number of $\mathrm{SL}(2,\mathbb{Z})$-equivalence classes of binary quadratic forms with fixed discriminant weighted by the order of their automorphism group (see Section \ref{sec 2}). Zagier in 1975 discovered the transformation properties of the generating function of Hurwitz class number \cite{DZ1}. Remarkably enough, his results show that in the modern framework of Zwegers, this generating function is a mock modular form of weight $3/2.$ Ramanujan's mock theta functions and Zagier's weight 3/2 mock modular form are some of the earliest examples of the theory. 
While there are excellent references \cite{HZ,Ono1} which provide an outline and discuss certain aspects of the proof  of this result, to our knowledge, a detailed step by step proof has not appeared in the literature yet. The aim of this paper is to provide such a proof. Even though the ideas are not new, we feel that this paper will be of benefit to those interested in understanding the proof of this important result in depth. \\\\
The article is organised as follows: in Section \ref{sec 2}, we recall the definitions of modular forms and mock modular forms and state Zagier's theorem. In Section \ref{sec 3} we review Cohen-Eisenstein series which are used to construct Zagier's mock modular forms and give detailed proofs of the modularity transformations and the analytic continuation. Finally, we prove Zagier's theorem (Theorem \ref{thm 1.3}) in Section \ref{sec 4}. 
\section{Preliminaries}
\label{sec 2}
Let $\mathbb{H}\coloneqq\{\tau=u+iv~|~u,v\in\mathbb{R}, v>0\}$ denote the complex upper half plane and $\mathrm{SL}(2,\mathbb{Z})$ be the group of $2\times 2$ matrices with integer entries and determinant 1. The group $\mathrm{SL}(2,\mathbb{Z})$ acts on $\mathbb{H}$ in the following way: for $\tau\in\mathbb{H},~~\gamma=\left(\begin{smallmatrix}
a&b\\c&d
\end{smallmatrix}\right)\in\mathrm{SL}(2,\mathbb{Z})$, we have 
\[
\gamma \tau=\frac{a\tau+b}{c\tau+d}.
\]
We usually restrict to the action of the subgroup $\Gamma_0(N)$ of $\mathrm{SL}(2,\mathbb{Z})$ called the Hecke congruence subgroup of level $N$  defined by 
\begin{equation*}
\Gamma_0(N)\coloneqq \left\{\begin{pmatrix}
a&b\\c&d\\
\end{pmatrix}\in \mathrm{SL}(2,\mathbb{Z}):c~\equiv ~0~(\text{mod} ~N)\right\}.
\end{equation*}
We now define modular forms and refer the reader to \cite{cohen} for details. 
\begin{defn}
A holomorphic function $f:\mathbb{H}\longrightarrow\mathbb{C}$ is called a holomorphic modular form (resp. weakly holomorphic modular form resp. cusp form) of weight $k\in\mathbb{Z}$  and level $N$ if $f(\gamma\tau)=(c\tau+d)^kf(\tau)$ for every $\gamma=\left(\begin{smallmatrix}
a&b\\c&d
\end{smallmatrix}\right)\in\Gamma_0(N)$ and $f$ is holomorphic (resp. meromorphic resp. vanishes) at all cusps of $\Gamma_0(N).$ 
\end{defn}
We have the following notation for the corresponding $\mathbb{C}-$vector spaces :
\[
\begin{split}
&M_k(\Gamma_0(N)) ~:~\text{space of holomorphic modular forms}\\
&M_k^!(\Gamma_0(N)) ~:~\text{space of weakly holomorphic modular forms}\\
&S_k(\Gamma_0(N)) ~:~\text{space of cusp forms.}
\end{split}
\] 
When $k\in\frac{1}{2}+\mathbb{Z},$ one defines modular forms with a slightly different transformation property (see Definition 1.2 \ref{def 1.1 (i)}). The definition of modular forms can be further generalised by relaxing the holomorphicity condition and requiring that our function be an eigenfunction of a certain differential operator. To this end, define the weight-$k$ $(\in\mathbb{R})$ hyperbolic Laplacian 
\begin{equation*}
\Delta_k\coloneqq -v^2\left(\frac{\partial^2}{\partial u^2}+\frac{\partial^2}{\partial v^2}\right)+ikv\left(\frac{\partial}{\partial u}+i\frac{\partial}{\partial v}\right)=-4v^2\frac{\partial}{\partial \tau}\frac{\partial}{\partial \bar\tau}+2ikv\frac{\partial}{\partial \bar\tau}.
\end{equation*}
Classically one begins by defining the classical Maass forms which are eigenfunctions of the weight $0$ hyperbolic Laplacian with eigenvalue of the form $\lambda=r(1-r)$ where $r$ is a parameter, and satisfy $f(\gamma\tau)=f$ for every $\gamma\in\Gamma_0(N)$. The main objects of this exposition are harmonic Maass forms and mock modular forms which we now define. We recall basic facts about harmonic Maass forms and refer the reader to  \cite{BF,Ono1} for details.
\theoremstyle{defn}
\begin{defn}
\normalfont
Let $k\in\frac{1}{2}\mathbb{Z}$. A smooth function (in real sense) $f:\mathbb{H}\rightarrow\mathbb{C}$ is called a weight-$k$ harmonic Maass form on $\Gamma_0(N)$ ($4|N$ if $k\in\frac{1}{2}\mathbb{Z}\setminus\mathbb{Z}$) if 
\begin{enumerate}[(i)]
\item For all $\begin{pmatrix}
a&b\\c&d\\
\end{pmatrix}\in \Gamma_0(N)$ and $\tau\in\mathbb{H}$, we have
\[
  f\left(\frac{az+b}{cz+d}\right) =
  \begin{cases}
                                   (cz+d)^kf(\tau) & \text{if $k\in \mathbb{Z}$} \\
                                   (\frac{c}{d})\varepsilon_d^{-2k}(cz+d)^kf(\tau) & \text{if $k\in \frac{1}{2}+\mathbb{Z}.$} \\  
  \end{cases}
\] 
Here $\left(\frac{c}{d}\right)$ is the Jacobi symbol and $\varepsilon_d=\sqrt{\left(\frac{-1}{d}\right)}$.
\label{def 1.1 (i)}
\item $\Delta_k(f)=0$.
\item There exists a polynomial $P_f(\tau)\in\mathbb{C}[q^{-1}]$ such that $f(\tau)-P_f(\tau)=O(e^{-\epsilon v})$ as $v\rightarrow\infty$ for some $\epsilon>0$. Similar conditions hold at other cusps.
\end{enumerate}
Here $\sqrt{\cdot}$ is the principal branch of square root mapping upper half plane to first quadrant and lower half plane to fourth quadrant.
If the third condition in the above Definition is replaced by $f(\tau)=O(e^{\epsilon v})$ for some $\epsilon>0$, then $f$ is said to be a harmonic Maass form of \emph{manageable growth}. Space of harmonic Maass forms of weight $k$ is denoted by $H_k(\Gamma_0(N))$ and that of harmonic Maass forms of manageable growth is denoted by $H_k^{!}(\Gamma_0(N))$. 
\end{defn} 
Any $f\in H_k^{!}(\Gamma_0(N))$ has a Fourier expansion \cite{Ono1} of the form  
\begin{equation}
f(\tau)=f(u+iv)=\sum\limits_{n>> -\infty}c_f^{+}(n)q^n+ c_f^{-}(0)v^{1-k}+\sum\limits_{\substack{n<< \infty\\n\neq 0}}c_f^{-}(n)\Gamma(1-k,-4\pi nv)q^n,~~q=e^{2\pi i\tau},
\label{1}
\end{equation}
where $\Gamma(s,z)$ is the incomplete gamma function defined as 
\begin{equation}
\Gamma(s,z)\coloneqq\int\limits_z^{\infty}e^{-t}t^s\frac{dt}{t}.
\label{2}
\end{equation}
The notation $\sum\limits_{n>> -\infty}$ means $\sum\limits_{n=\alpha_f}^{\infty}$ for some $\alpha_f\in\mathbb{Z}$ and $\sum\limits_{n<<\infty}$ is defined similarly. 
We call 
\begin{equation*}
f^+(\tau)=\sum\limits_{n>> -\infty}c_f^{+}(n)q^n,
\end{equation*}
the \emph{holomorphic part} of $f$ and 
\begin{equation*}
f^-(\tau)=c_f^{-}(0)v^{1-k}+\sum\limits_{\substack{n<< \infty\\n\neq 0}}c_f^{-}(n)\Gamma(1-k,-4\pi nv)q^n,
\end{equation*}
 the \emph{nonholomorphic part} of $f$.
Define the shadow operator as 
\begin{equation*}
\xi_k\coloneqq 2iv^k\compconj{\frac{\partial}{\partial\bar\tau}}.
\end{equation*} 
It is related to the Laplacian operator by 
\begin{equation}
\Delta_k=-\xi_{2-k}\circ\xi_{k}.
\label{3}
\end{equation}
It is known that for $f\in H_{2-k}^{!}(\Gamma_0(N))$, we have $\xi_{2-k}(f)\in M_{k}^!(\Gamma_0(N))$. Furthermore, given the Fourier expansion of $f$ as in \eqref{1}, we have 
\begin{equation}
\xi_{2-k}(f(\tau))=\xi_{2-k}(f^-(\tau))=(k-1)\compconj{c_f^-(0)}-(4\pi)^{k-1}\sum\limits_{n>>-\infty}\compconj{c_f^-(-n)}n^{k-1}q^n,
\label{4}
\end{equation}
where $\compconj{c_f^-(0)}$ denoted the complex conjugate of $c_f^-(0)$. The image of $f$ under $\xi$ is called the shadow of $f$. Moreover this map is surjective. 
\begin{defn}
A \emph{mock modular form} of weight $(2-k)$ is the holomorphic part $f^+$ of a harmonic Maass form of weight $(2-k)$ for which $f^-$ is non trivial. The weakly holomorphic modular form $\xi_{2-k}(f)$ is called the shadow of the mock modular form $f^+$. 
\end{defn}
The first example of mock modular form comes from the weight-2 nonholomorphic Eisenstein series defined as 
\[
E_2^*(\tau)=E_2(\tau)-\frac{3}{\pi v},
\]
where 
\[
E_{2}(\tau)=1-24 \sum_{n=1}^{\infty} \sigma_{1}(n) q^{n},
\]
with $\sigma_k$ being the $k$th divisor sum. One can easily check that $E_2^*$ transforms as modular form and has proper growth at the cusps \cite{Ono1}. Thus $E_2^*\in H_2^!(\mathrm{SL}(2,\mathbb{Z}))$ and hence $E_2$ is a mock modular form of weight $2$ for $\mathrm{SL}(2,\mathbb{Z}).$  
\\\\
As mentioned in the introduction, generating functions of certain generalisations of class numbers turn out to be modular forms as we explain below. We define the generalised Cohen class number as follows: 
Let $r,N\in\mathbb{N}$. Following \cite{Co}, we define
$$
h(r,N)=\left\{\begin{array}{l}
(-1)^{[r / 2]}(r-1) ! N^{r-1 / 2} 2^{1-r} \pi^{-r} L\left(r, \chi_{(-1) r N}\right) \text { if }(-1)^{r} N \equiv 0 \text { or } 1(\bmod 4) \\
0 \text { if }(-1)^{r} N \equiv 2 \text { or } 3(\bmod 4),
\end{array}\right.
$$
where we write $\chi_{D}$ for the character $\chi_{D}(d)=\left(\frac{D}{d}\right)$. Using $h(r,N)$, we define the Cohen class number 
$$
H(r, N)=\left\{\begin{array}{lll}
\sum\limits_{d^{2} \mid N} h\left(r, N / d^{2}\right) & \text { if }(-1)^{r} N \equiv 0 & \text { or } \quad 1(\bmod~4), \quad N>0 \\
\zeta(1-2 r) \quad \text { if } & N=0 \\
0 & \text { otherwise }.
\end{array}\right.
$$
As a consequence of Dirichlet's class number formula (cf. remark (c) in Section 2 of \cite{Co}) we have 
\begin{equation}
 H(r,n)  =
  \begin{cases}
                                   0 & \text{if $(-1)^rn~\equiv~2,3(\text{mod}~4)$} \\
                                   L(1-r,\chi_d)T^{\chi_d}_r(f)& \text{if $(-1)^rn~\equiv~0,1(\text{mod}~4), n>0$} \\
\zeta(1-2r)& \text{if $n=0$}, \\  
  \end{cases}
\label{5} 
\end{equation}
where $L(s,\chi)$ is the Dirichlet $L$-function, $\chi_d$ is the character associated to the quadratic field $\mathbb{Q}(\sqrt{(-1)^rn})$ and $(-1)^rn=df^2$ and is defined using the Jacobi symbol as follows:
\begin{equation}
\chi_d(a)=\left(\frac{d}{a}\right),
\label{6}
\end{equation} 
and $T^{\chi_d}_s(f)$ is defined as follows:
\begin{equation}
T^{\chi_d}_s(f)\coloneqq \sum\limits_{a|f}\mu(a)\chi_d(a)a^{k-1}\sigma_{2s-1}(f/a),
\label{7} 
\end{equation}
where $\mu(a)$ is M\"{o}bius function and $\sigma_{k}$ is the divisor function of order $k$. With this definition, we consider the function 
\begin{equation}
H_{k+1/2}(\tau)=\sum_{n=0}^{\infty}H(k,n)q^n,~~k\geq 2.
\label{eq 8}
\end{equation}
In \cite{Co}, Cohen proved that $H_{k+1/2}(\tau)\in M_{k+1/2}(\Gamma_0(4)).$ 
\theoremstyle{thm}
\begin{defn}(Hurwitz class number)
For integer $N\geq 0$, the Hurwitz class number  is defined as follows : $H(0)=-\frac{1}{12}$. If $N\equiv 1~\text{or}~2(\text{mod}~4)$ then $H(N)=0$. Otherwise $H(N)$ is the number of classes of binary quadratic forms of discriminant $-N$, except that those classes which have a representative which is a multiple of $x^2+y^2$ should be counted with weight $1/2$ and those which have a representative which is a multiple of $x^2+xy+y^2$ should be counted with weight $1/3$. 
\end{defn}
\begin{remark}
\label{remark 2.5}
As a result of Dirichlet's class number formula, for $k=1$, $H(k,n)$ reduces to the Hurwitz class number.
\end{remark} 
Zagier's results show that  $H_{3/2}(\tau)$ as defined in \eqref{eq 8} is not a modular form, rather its transformation properties can be understood in terms of mock modular forms. The precise statement is formulated (cf. \cite[Theorem 6.3]{Ono1}) in the theorem below. 
\theoremstyle{thm}
\begin{thm}\emph{(Zagier, \cite{DZ1})}
Let $H(n)$ be the Hurwitz class number. Then the function
\begin{equation*}
\mathcal{H}(\tau)\coloneqq -\frac{1}{12}+\sum\limits_{n=1}^{\infty}H(n)q^n+\frac{1}{4\sqrt{\pi}}\sum\limits_{n=1}^{\infty}n\Gamma\left(-\frac{1}{2},4\pi n^2v\right)q^{-n^2}+\frac{1}{8\pi\sqrt{v}},
\end{equation*}
where $\tau=u+iv$ is a $3/2-$weight harmonic Maass form of manageable growth on $\Gamma_0(4)$. Moreover we have $\xi_{3/2}(\mathcal{H})=-\frac{1}{16}\Theta$ where $\Theta(\tau)\coloneqq\sum\limits_{n\in\mathbb{Z}}q^{n^2}$ is the classical theta function.
\label{thm 1.3}
\end{thm}
This theorem says that the generating function of Hurwitz class numbers is a mock modular forms of weight $3/2$ for $\Gamma_0(4).$ We will now give a detailed proof of this theorem.  
\section{Cohen-Eisenstein Series}
\label{sec 3}
Let $k\in\mathbb{Z}$ and $\tau=u+iv\in\mathbb{H}$. Define the weight-$\left(k+\frac{1}{2}\right)$ Eisenstein series as
\begin{equation}
\begin{split}
& E_{k+1/2}(\tau)\coloneqq \sum_{\substack{m=1\\m~\text{odd}}}^{\infty}\sum_{\substack{n=-\infty\\\text{gcd}(n,m)=1}}^{\infty}\frac{\left(\frac{n}{m}\right)\varepsilon_m^{-2k-1}}{(m\tau+n)^{k+1/2}},\quad\text{and}\\
& F_{k+1/2}(\tau)\coloneqq \tau^{-k-1/2}E_{k+1/2}\left(\frac{-1}{4\tau}\right).
\end{split}
\end{equation}
The series converges absolutely for $k>3/2$. Shimura proved that $E_{k+1/2},F_{k+1/2}$ are modular forms of weight $\left(k+\frac{1}{2}\right)$ for $\Gamma_0(4)$ and computed their  Fourier expansion \cite{Sh1,Sh2}.
Following Cohen \cite{Co}, we define the linear combination of these two Eisenstein series 
\begin{equation}
H_{k+1/2}(\tau)\coloneqq \frac{\zeta(1-2k)}{2^{2k+1}}\left[(1+i^{2k+1})E_{k+1/2}(\tau)+i^{2k+1}F_{k+1/2}(\tau)\right].
\label{9}
\end{equation}
For $k\geq 2$, the two expressions for $H_{k+1/2}$ as in \eqref{eq 8} and \eqref{9} coincide. Zagier proved that the function $\mathcal{H}(\tau)=H_{3/2}(\tau)$ where $H_{3/2}(\tau)$ is as in \eqref{9}. To prove the transformation property of $H_{3/2}$, one defines the nonholomorphic Eisenstein series as follows: for $\tau\in\mathbb{H}, s\in\mathbb{C}$, put
\begin{equation}
\begin{split}
& E_{k+1/2,s}(\tau)\coloneqq \sum_{\substack{m=1\\m~\text{odd}}}^{\infty}\sum_{\substack{n=-\infty\\\text{gcd}(n,m)=1}}^{\infty}\frac{\left(\frac{n}{m}\right)\varepsilon_m^{-2k-1}}{(m\tau+n)^{k+1/2}|m\tau+n|^{2s}},\\
& F_{k+1/2,s}(\tau)\coloneqq \tau^{-k-1/2}|\tau|^{2s}E_{k+1/2,s}\left(\frac{-1}{4\tau}\right),\\&
H_{k+1/2,s}(\tau)\coloneqq \frac{\zeta(1-2k)}{2^{2k+1}}\left[(1+i^{2k+1})E_{k+1/2,s}(\tau)+i^{2k+1}F_{k+1/2,s}(\tau)\right].
\end{split}
\end{equation}
Each series converges absolutely for Re$(s)>\frac{3}{4}-\frac{k}{2}$. This restriction is assumed wherever we use the series representation of $E_{k+1/2,s},F_{k+1/2,s}$ and $H_{k+1/2,s}$. In subsequent sections, we will frequently use absolute convergence to rearrange the series as required.
\subsection{Modularity transformation for $H_{k+1/2,s}$}
We will show that for $\begin{pmatrix}
a&b\\c&d
\end{pmatrix}\in \Gamma_0(4)$,
\begin{equation}
\begin{split}
E_{k+1/2,s}\left(\frac{a\tau+b}{c\tau+d}\right)=\left(\frac{c}{d}\right)\varepsilon_m^{-2k-1}(c\tau+d)^{k+1/2}|c\tau+d|^{2s}E_{k+1/2,s},\quad\text{and}\\F_{k+1/2,s}\left(\frac{a\tau+b}{c\tau+d}\right)=\left(\frac{c}{d}\right)\varepsilon_m^{-2k-1}(c\tau+d)^{k+1/2}|c\tau+d|^{2s}F_{k+1/2,s}.
\end{split}
\label{11}
\end{equation}
Thus $H_{k+1/2,s}$ also satisfies similar modularity transformation with respect to $\Gamma_0(4)$. 
Let us introduce the automorphy factor of the theta function:
\begin{equation*}
J(\gamma,\tau)\coloneqq v_{\vartheta}(\gamma)j(\gamma,\tau)^{1/2},
\end{equation*}
where $v_{\vartheta}(\gamma)$ is the multiplier system for theta function. The theta multiplier $v_{\vartheta}(\gamma)$ and the factor $j(\gamma,\tau)$ are defined as 
\begin{equation*}
v_{\vartheta}(\gamma)\coloneqq \left(\frac{c}{d}\right)\varepsilon_d^{-1}~~\text{and}~~j(\gamma,\tau)\coloneqq c\tau+d,~~~~\gamma=\begin{pmatrix}
a&b\\c&d
\end{pmatrix}\in \Gamma_0(4),~\tau\in\mathbb{H}.
\end{equation*}
Since $J(\gamma,\tau)$ satisfies the cocycle condition: 
\begin{equation*}
J(\gamma_1\gamma_2,\tau)=J(\gamma_1,\gamma_2\tau)J(\gamma_2,\tau),
\end{equation*}
the multiplier system must satisfy 
\begin{equation}
v_{\vartheta}(\gamma_1\gamma_2)=\Sigma(\gamma_1,\gamma_2)v_{\vartheta}(\gamma_1)v_{\vartheta}(\gamma_2),
\label{eq:multcocycle}
\end{equation}
where 
\begin{equation*}
\Sigma(\gamma_1,\gamma_2)=j(\gamma_1,\gamma_2\tau)^{1/2}j(\gamma_2,\tau)^{1/2}j(\gamma_1\gamma_2,\tau)^{-1/2}\in\{\pm 1\}.
\end{equation*}
is independent of $\tau\in\mathbb{H}$.
Let us recall some properties of the Jacobi symbol which will be useful in our computations.
\theoremstyle{thm}
\begin{thm}\emph{(Proposition 3.4.7, \cite{cohen})}
For $a,m,n\in\mathbb{Z}$ we have \\
$(a)$ $(Vanishing)$
\[
\left(\frac{a}{n}\right)=0 \quad \text { if and only if } \quad \operatorname{gcd}(a, n)>1;
\]
$(b)$ $(Multiplicativity)$
\[
\left(\frac{a}{m}\right)\left(\frac{a}{n}\right)=\left(\frac{a}{m n}\right) \quad \text { and } \quad\left(\frac{a}{n}\right)\left(\frac{c}{n}\right)=\left(\frac{a c}{n}\right);
\]
$(c)$ $(\text {Periodicity})$ If $n>0$ or $n<0$ and $\operatorname{sign}(a)=\operatorname{sign}(c),$ then
\[
\left(\frac{a}{n}\right)=\left(\frac{c}{n}\right) \quad \text { if } a \equiv c \bmod \left\{\begin{array}{ll}
4 n & \text { if } n \equiv 2(\bmod 4) \\
n & \text { otherwise; }
\end{array}\right.
\]
Furthermore, if $a \neq 0$ and $a \neq 3(\bmod 4),$ then
\[
\left(\frac{a}{n}\right)=\left(\frac{a}{m}\right) \quad \text { if } n \equiv m \bmod \left\{\begin{array}{ll}
4|a| & \text { if } a \equiv 2(\bmod 4) \\
a & \text { if } a \equiv 0,1(\bmod 4).
\end{array}\right.
\]
\label{thm 2.1}
\end{thm}
We will prove the modularity of $E_{k+1/2,s}$ in a series of Lemmas.
\begin{lemma}
\begin{enumerate}[(i)]
\item For $z\in\mathbb{C}\setminus\mathbb{R}$ and $k\in\frac{1}{2}\mathbb{Z}$, we have
\[
(-z)^{k}=z^{k}(-i~\mathrm{sgn}(\mathrm{Im}(z)))^{2k}.
\] 
\item For any $z\in\mathbb{H}$ and $k\in\mathbb{Z}+\frac{1}{2}$, we have 
\[
\begin{split}
z^{-k}(az+b)^{k}=\begin{cases}
\left(a+\frac{b}{z}\right)^{k}&a>0\text{ or }b>0\\i^{4k}\left(a+\frac{b}{z}\right)^{k}&a<0\text{ and }b<0.
\end{cases}
\end{split}
\]
In particular, if $a<0$, then $$z^{-k}(az+b)^{k}=\mathrm{sign}(b)\left(a+\frac{b}{z}\right)^{k},$$ and when $b<0$, then $$z^{-k}(az+b)^{k}=\mathrm{sign}(a)\left(a+\frac{b}{z}\right)^{k}.$$
\item For any matrix $g=\begin{pmatrix}
a&b\\c&d
\end{pmatrix} \in \Gamma_0(4)$ with $a<0$ and $c>0$, and $\frac{5}{2}\leq k\in\frac{1}{2}\mathbb{Z}\setminus\mathbb{Z}$ we have that
\begin{equation*}
-\mathrm{sign}(b)\left(\frac{-b}{a}\right)\varepsilon_a^{-2k}=\left(\frac{c}{d}\right)\varepsilon_d^{-2k}.
\end{equation*} 
\end{enumerate}
\label{lemma:negsqrootcomp}
\begin{proof}
(i) First note that when $k\in\mathbb{Z},$ then the statement is easy to prove. When $z$ is real then also the statement is easily checked to be true. Suppose $z\in\mathbb{H},$ then
\[
\begin{aligned}
(-z)^{k} &=e^{k\ln (-z)} \\
&=e^{k[\ln |z|+i(-\pi+\theta)]} \\
&=e^{k \ln |z|} \cdot e^{ik(-\pi+\theta)} \\
&=e^{k(\ln |z|+i \theta)} e^{-i \pi k} \\
&=z^{k}(-i)^{2k}.
\end{aligned}
\] 
Next if $z\in\mathbb{C}\backslash\mathbb{H}$, then similar calculation gives $(-z)^k=i^{2k}z^k.$\\\\
(ii) If $a>0$, then $az+b\in\mathbb{H}$ and $z^{-1}$ is in the lower-half plane, thus in the principal branch of squareroot, we can combine $z^{-k}(az+b)^{k}$ to get
\[
z^{-k}(az+b)^{k}=\left(\frac{az+b}{z}\right)^{k}=\left(a+\frac{b}{z}\right)^{k}.
\]
When $b>0$, then we see that $a+b/z$ is in the lower-half plane and $z\in\mathbb{H}$, hence we can combine $z^{k}(a+b/z)^{k}$ to get 
\[
z^k\left(a+\frac{b}{z}\right)^{k}=(az+b)^k.
\]
Finally when both $a$ and $b$ is negative, then using the first part of this lemma, we have 
\[
\begin{split}
z^{-k}(az+b)^{k}&=z^{-k}(-(-az-b))^k=(-i)^{2k}\left(\frac{1}{z}\right)^{k}(-az-b)^{k}\\&=(-i)^{2k}\left(-\left(a+\frac{b}{z}\right)\right)^k=(-i)^{4k}\left(a+\frac{b}{z}\right)^{k},
\end{split}
\]
where we used the fact that for both $a$ and $b$ is negative, $-az-b$ and $a+b/z$ are in lower-half plane. The second statement is now obvious.
\\\\
(iii) The proof follows by considering the modularity of the Eisenstein series $\widetilde{F}_{k}\coloneqq (2z)^{-k}\widetilde{E}_{k}\left(-\frac{1}{4z}\right)$ for $\frac{5}{2}\leq k\in\frac{1}{2}\mathbb{Z}\setminus\mathbb{Z}$ where $\widetilde{E}_{k}$ defined as: 
\begin{equation*}
\widetilde{E}_{k}(z)=\sum_{\gamma\in\Gamma_{\infty}\setminus\Gamma_0(4)}J(\gamma,z)^{-2k}.
\end{equation*} 
For a proof of modularity of $\widetilde{E}_{k}$ and $\widetilde{F}_{k}$, we refer the reader to $\S$2 of Chapter IV of \cite{Koblitz}. We have $\widetilde{F}_{k}(gz)=J(g,z)^{2k}\widetilde{F}_{k}(z)$. Using the definition of $\widetilde{E}_{k}$, we get
\begin{equation*}
\widetilde{F}_{k}(gz)=2^{-k}\left(\frac{az+b}{cz+d}\right)^{-k}\widetilde{E}_{k}\left(-\frac{1}{4\left(\frac{az+b}{cz+d}\right)}\right)=2^{-k}\left(\frac{az+b}{cz+d}\right)^{-k}\widetilde{E}_{k}\left(g'\left(-\frac{1}{4z}\right)\right),
\end{equation*} 
where $g'=\begin{pmatrix}
d&-c/4\\-4b&a
\end{pmatrix}$. Since we have 
$\widetilde{E}_{k}\left(g'\left(-\frac{1}{4z}\right)\right)=J\left(g',-\frac{1}{4z}\right)^{2k}\widetilde{E}_{k}\left(-\frac{1}{4z}\right),$ we get 
\begin{equation*}
J(g,z)^{2k}z^{-k}=J\left(g',-\frac{1}{4z}\right)^{2k}\left(\frac{az+b}{cz+d}\right)^{-k}.
\end{equation*}
Substituting the expression for $J$, we get
\begin{equation*}
\begin{split}
\left(\frac{c}{d}\right)\varepsilon_d^{-2k}(cz+d)^{k}z^{-k}&=\left(\frac{-4b}{a}\right)\varepsilon_a^{-2k}\left(\frac{b}{z}+a\right)^{k}\left(\frac{az+b}{cz+d}\right)^{-k}\\&=\left(\frac{-b}{a}\right)\varepsilon_a^{-2k}\mathrm{sign}(b)z^{-k}\left(az+b\right)^{k}\left(\frac{az+b}{cz+d}\right)^{-k},
\end{split}
\end{equation*}
where we used Lemma \ref{lemma:negsqrootcomp} (ii). Next we note that $(az+b)^{-1},(cz+d)\in \mathbb{H}$. Thus we have 
\[
\left(-\frac{1}{(-az-b)}\right)^{k}(cz+d)^k=i^{2k}\left(\frac{1}{(-az-b)}\right)^{k}(cz+d)^k=i^{2k}\left(-\frac{cz+d}{az+b}\right)^k=i^{4k}\left(\frac{cz+d}{az+b}\right)^k,
\]
where we used Lemma \ref{lemma:negsqrootcomp} (i) along with the fact that $(-az-b)^{-1}$ and $(cz+d)/(az+b)$ are in the lower half plane. 
This gives 
 \begin{equation*}
\begin{split}
\left(\frac{c}{d}\right)\varepsilon_d^{-2k}(cz+d)^{k}z^{-k}&=-\mathrm{sign}(b)\left(\frac{-b}{a}\right)\varepsilon_a^{-2k}z^{-k}(cz+d)^{k}.
\end{split}
\end{equation*}
Since this holds for every $z\in\mathbb{H}$, we get our result.
\end{proof}
\end{lemma}
\theoremstyle{thm}
\begin{lemma}
For any $\gamma,g,S\in\mathrm{SL}(2,\mathbb{R})$, we have that 
\begin{equation*}
\Sigma(S^{-1},g)\Sigma(S^{-1}g,\gamma)=\Sigma(S^{-1},g\gamma)\Sigma(g,\gamma).
\end{equation*}
In particular, if $g=S\eta S^{-1}$ with $\eta=\left(\begin{smallmatrix}1&n\\0&1\end{smallmatrix}\right)$ for some $n\in\mathbb{Z}$,  then 
\[
\Sigma(S^{-1},\gamma)=\Sigma(S^{-1},g\gamma)\Sigma(g,\gamma).
\] 
\label{lemma 2.3} 
\begin{proof}
By definition, we have 
\[
\begin{split}
\Sigma(S^{-1},g)\Sigma(S^{-1}g,\gamma)&=j(S^{-1},g\tau)^{1/2}j(g,\tau)^{1/2}j(S^{-1}g,\tau)^{-1/2}\\&\hspace{3cm}\times j(S^{-1}g,\gamma\tau)^{1/2}j(\gamma,\tau)^{1/2}j(S^{-1}g\gamma,\tau)^{-1/2}.
\end{split}
\]
and 
\[
\begin{split}
\Sigma(S^{-1},g\gamma)\Sigma(g,\gamma)&=j(S^{-1},g\gamma\tau)^{1/2}j(g\gamma,\tau)^{1/2}j(S^{-1}g\gamma,\tau)^{-1/2}\\&\hspace{3cm}\times j(g,\gamma\tau)^{1/2}j(\gamma,\tau)^{1/2}j(g\gamma,\tau)^{-1/2}\\&=\Sigma(S^{-1}g,\gamma)j(S^{-1}g,(\gamma\tau))^{-1/2}j(S^{-1},g(\gamma\tau))^{1/2}j(g,(\gamma\tau))^{1/2}\\&=\Sigma(S^{-1},g)\Sigma(S^{-1}g,\gamma),
\end{split}
\]
where in the last step we used the fact that $\Sigma(\gamma_1,\gamma_2)$ is independent of $\tau\in\mathbb{H}$. To prove the second relation, it suffices to prove that $\Sigma(S^{-1},g)\Sigma(\eta S^{-1},\gamma)=\Sigma(S^{-1},\gamma)$. Indeed, by definition we have
\[
\begin{split}
\Sigma(S^{-1},S\eta S^{-1})&=j(S^{-1},S\eta S^{-1}\tau)^{1/2}j(S\eta S^{-1},\tau)^{1/2}j(\eta S^{-1},\tau)^{-1/2}\\&=\sqrt{-c\left(\frac{(1-acn)\tau+a^2n}{-c^2n\tau+1+acn}+a\right)}\frac{\sqrt{-c^2n\tau+1+acn}}{\sqrt{-c\tau+a}},
\end{split}
\]
where 
\[
S=\begin{pmatrix}
a&b\\c&d
\end{pmatrix}\quad\text{so that}\quad S^{-1}=\begin{pmatrix}
d&-b\\-c&a
\end{pmatrix},\quad S^{-1}\eta S=\begin{pmatrix}
1-acn & a^2n\\-c^2n & 1+acn
\end{pmatrix}.
\]
Simplifying and using Lemma \ref{lemma:negsqrootcomp} (ii), we get 
\[
\begin{split}
\Sigma(S^{-1},S\eta S^{-1})&=\sqrt{\frac{-c\tau+a}{cn(-c\tau+a)+1}}\frac{\sqrt{cn(-c\tau+a)+1}}{\sqrt{-c\tau+a}}\\&=\frac{\sqrt{-c\tau+a}}{\sqrt{cn(-c\tau+a)+1}}\frac{\sqrt{cn(-c\tau+a)+1}}{\sqrt{-c\tau+a}}\\&=1.
\end{split}
\]
Finally since left multiplication by $\eta$ to a matrix does not change the bottom row of the matrix, we have
\[
\Sigma(\eta S^{-1},\gamma)=\Sigma(S^{-1},\gamma).
\] 
\end{proof}
\end{lemma}

\theoremstyle{thm}
\begin{lemma}
For $k\in\mathbb{Z}$ and Re$(s)>\frac{3}{4}-\frac{k}{2}$, we have 
\begin{equation*}
E_{k+1/2,s}(\tau)=i^{-2k-1}\sum_{\gamma\in \Gamma_0(4)_0\setminus \Gamma_0(4)}\compconj{v_{\vartheta}(\gamma)}^{2k+1}\Sigma(S^{-1},\gamma)j(S^{-1}\gamma,\tau)^{-k-1/2}|j(S^{-1}\gamma,\tau)|^{-2s},
\end{equation*}
where $\Gamma_0(4)_0$ is the stabiliser in $\Gamma_0(4)$ of the cusp 0 of $\Gamma_0(4)$ and $S=\begin{pmatrix}
0&-1\\1&0
\end{pmatrix}$.
\label{lemma 2.4} 
\begin{proof}
We first show that the sum on the right hand side is well defined, that is the summands $\compconj{v_{\vartheta}(\gamma)}^{2k+1}\Sigma(S^{-1},\gamma)j(S^{-1}\gamma,\tau)^{-k-1/2}|j(S^{-1}\gamma,\tau)|^{-2s}$ are independent of the choice of representatives. Observe that $\Gamma_0(4)_0=S\Gamma_{\infty}S^{-1}$ where 
\begin{equation*}
\Gamma_{\infty}\coloneqq \left\{\pm\begin{pmatrix}
1&n\\0&1
\end{pmatrix}~:~n\in\mathbb{Z}\right\}
\end{equation*} 
is the stabiliser of $i\infty$. This is because $S(i\infty)=0.$  Let $\gamma'=\beta\gamma$ where $\beta\in\Gamma_0(4)_0$. Then $\beta=S\eta S^{-1}$ for some $\eta\in\Gamma_{\infty}$. We need to show that 
\[
\begin{split}
\compconj{v_{\vartheta}(\gamma)}^{2k+1}&\Sigma(S^{-1},\gamma)j(S^{-1}\gamma,\tau)^{-k-1/2}|j(S^{-1}\gamma,\tau)|^{-2s}\\&=\compconj{v_{\vartheta}(\gamma')}^{2k+1}\Sigma(S^{-1},\gamma')j(S^{-1}\gamma',\tau)^{-k-1/2}|j(S^{-1}\gamma',\tau)|^{-2s}.
\end{split}
\]
Indeed, since left multiplication by $\eta$ to a matrix does not change the bottom row of the matrix, we have 
\[
j(S^{-1}\gamma',\tau)=j(\eta S^{-1}\gamma,\tau)=j(S^{-1}\gamma,\tau).
\]
Next using \eqref{eq:multcocycle} and Lemma \ref{lemma 2.3}, we have
\[
\begin{split}
& v_{\vartheta}(\gamma')=v_{\vartheta}(\beta\gamma)=\Sigma(\beta,\gamma)v_{\vartheta}(\beta)v_{\vartheta}(\gamma)\\&\Sigma(\beta,\gamma)\Sigma(S^{-1},\gamma')=\Sigma(\beta,\gamma)\Sigma(S^{-1},\beta\gamma)=\Sigma(S^{-1},\gamma).
\end{split}
\]
It now follows that the sum is well defined from the fact that $\Sigma(\beta,\gamma)$ is real. Let us now simplify the right hand side. Now in the sum on the right hand side, substitute $S^{-1}\gamma=g$. We get
\begin{equation*}
\begin{split}
\sum_{\gamma\in\Gamma_0(4)_0\setminus \Gamma_0(4)}\compconj{v_{\vartheta}(\gamma)}^{2k+1}&\Sigma(S^{-1},\gamma)j(S^{-1}\gamma,\tau)^{-k-1/2}|j(S^{-1}\gamma,\tau)|^{-2s}\\&=\sum_{\gamma \in S\Gamma_{\infty}S^{-1}\setminus \Gamma_0(4)}\compconj{v_{\vartheta}(\gamma)}^{2k+1}\Sigma(S^{-1},\gamma)j(S^{-1}\gamma,\tau)^{-k-1/2}|j(S^{-1}\gamma,\tau)|^{-2s}\\ &=\sum_{g \in \Gamma_{\infty}\setminus S^{-1}\Gamma_0(4)}\compconj{v_{\vartheta}(Sg)}^{2k+1}\Sigma(S^{-1},Sg)j(g,\tau)^{-k-1/2}|j(g,\tau)|^{-2s}.
\end{split}
\end{equation*} 
To see this observe that if $\gamma_1=(S\gamma_{\infty}S^{-1})\gamma_2$ for some $\gamma_{\infty}\in\Gamma_{\infty}$, then $(S^{-1}\gamma_1)=\gamma_{\infty}(S^{-1}\gamma_2)$. Thus if we substitute $S^{-1}\gamma$ by $g$, then $g$ runs over the coset representatives of $\Gamma_{\infty}\setminus S^{-1}\Gamma_0(4).$ Now if 
\begin{equation*}
g=\begin{pmatrix}
a&b\\c&d
\end{pmatrix}\in S^{-1}\Gamma_0(4),~~\text{then}~~Sg=\begin{pmatrix}
-c&-d\\a&b
\end{pmatrix}\in\Gamma_0(4). 
\end{equation*}
Note that $c$ is odd and $d$ is an arbitrary integer. Moreover since 
\[
\begin{pmatrix}
1&n\\0&1
\end{pmatrix}\begin{pmatrix}
a&b\\c&d
\end{pmatrix}=\begin{pmatrix}
a+cn&b+dn\\c&d
\end{pmatrix},
\]
so we can choose $n\in\mathbb{Z}$ such that $a+cn>0$. Thus we choose representatives $g=\left(\begin{smallmatrix}a&b\\c&d\end{smallmatrix}\right)\in\Gamma_{\infty}\setminus S^{-1}\Gamma_0(4)$ such that $a,c>0$.  Now using Lemma \ref{lemma:negsqrootcomp} (iii) for $Sg=\left(\begin{smallmatrix}
-c&-d\\a&b
\end{smallmatrix}\right)\in\Gamma_0(4)$ and taking complex conjugate, we have 
\[
\begin{split}
v_{\vartheta}(Sg)=\left(\frac{a}{b}\right)\varepsilon_{b}^{2k+1}&=-\text{sign}(-d)\left(\frac{d}{-c}\right)\varepsilon_{-c}^{2k+1}\\&=\text{sign}(d)^2i^{2k+1}\left(\frac{d}{c}\right)\varepsilon_{c}^{-2k-1}\\&=i^{2k+1}\left(\frac{d}{c}\right)\varepsilon_{c}^{-2k-1},
\end{split}
\]
where we used the following facts:
\begin{enumerate}
\item Jacobi symbol is completely multiplicative in both entries.
\item $\left(\frac{d}{-1}\right)=\text{sign}(d)$.
\item $\varepsilon_{-c}=i\varepsilon_c^{-1}$.
\end{enumerate}  
Moreover there is a one-to-one correspondence between the coset representatives of $\Gamma_{\infty}\setminus S^{-1}\Gamma_0(4)$ and the set 
\begin{equation*}
W=\{(c,d)\in\mathbb{Z}^2:c>0~\text{odd}, d\in\mathbb{Z},~\text{gcd}(c,d)=1, c=1~\text{if}~d=0\}.
\end{equation*} 
The condition $c>0$ is to ensure that we only takes one pair out of $\pm(c,d)$ since $\Gamma_{\infty}$ identifies these two. Noting that $j(g,\tau)^{-k-1/2}=(c\tau+d)^{-k-1/2}$ for $g=\left(\begin{smallmatrix}
a&b\\c&d
\end{smallmatrix}\right)$, we get 
\begin{equation*}
\begin{split}
\sum_{\gamma\in\Gamma_0(4)_0\setminus \Gamma_0(4)}\compconj{v_{\vartheta}(\gamma)}^{2k+1}&\Sigma(S^{-1},\gamma)j(S^{-1}\gamma,\tau)^{-k-1/2}|j(S^{-1}\gamma,\tau)|^{-2s}\\&=i^{2k+1}\sum_{\substack{c=1\\c~\text{odd}}}^{\infty}\sum_{\substack{d\in\mathbb{Z}\\\text{gcd}(c,d)=1}}\left(\frac{d}{c}\right)\varepsilon_{c}^{-2k-1}(c\tau+d)^{-k-1/2}|c\tau+d|^{-2s}
\end{split}
\end{equation*}
where we used the fact that $\left(\frac{a}{b}\right)^{2k+1}=\left(\frac{a}{b}\right),\compconj{\varepsilon_{b}^{-1}}=\varepsilon_b$.
\end{proof}
\end{lemma}
\theoremstyle{thm}
\begin{thm}
For $\gamma\in\Gamma_0(4)$, we have 
\begin{equation*}
E_{k+1/2,s}(\gamma\tau)=v_{\vartheta}^{2k+1}j(\gamma,\tau)^{k+1/2}|j(\gamma,\tau)|^{2s}E_{k+1/2,s}(\tau).
\end{equation*}
\label{thm 2.5}
\begin{proof}
We have 
\begin{equation*}
E_{k+1/2,s}(\gamma\tau)=M\sum_{g\in\Gamma_0(4)_0\setminus \Gamma_0(4)}\compconj{v_{\vartheta}(g)}^{2k+1}\Sigma(S^{-1},g)j(S^{-1}g,\gamma\tau)^{-k-1/2}|j(S^{-1}g,\gamma\tau)|^{-2s},
\end{equation*}
where $M=(-1)^{-k+1/2}.$ Now note that \[
\begin{split}&|j(S^{-1}g,\gamma\tau)^{1/2}|=|j(S^{-1}g\gamma,\tau)^{1/2}||j(\gamma,\tau)^{-1/2}|\\&j(S^{-1}g,\gamma\tau)^{1/2}=\Sigma(S^{-1}g,\gamma)j(S^{-1}g\gamma,\tau)^{1/2}j(\gamma,\tau)^{-1/2}.
\end{split}
\] 
Using this and Lemma \ref{lemma 2.3}, we get 
\begin{equation*}
\begin{split}
& E_{k+1/2,s}(\gamma\tau)=j(\gamma,\tau)^{k+1/2}|j(\gamma,\tau)|^{2s}M\sum_{g\in\Gamma_0(4)_0\setminus \Gamma_0(4)}\compconj{v_{\vartheta}(g)}^{2k+1}\Sigma(S^{-1},g)\Sigma(S^{-1}g,\gamma)\\&\hspace{10cm}\times j(S^{-1}g\gamma,\tau)^{-k-1/2}|j(S^{-1}g\gamma,\tau)|^{-2s}\\&=j(\gamma,\tau)^{k+1/2}|j(\gamma,\tau)|^{2s}M\sum_{g\in\Gamma_0(4)_0\setminus \Gamma_0(4)}\compconj{v_{\vartheta}(g)}^{2k+1}\Sigma(g,\gamma)\Sigma(S^{-1},g\gamma)j(S^{-1}g\gamma,\tau)^{-k-1/2}\\&\hspace{10cm}\times|j(S^{-1}g\gamma,\tau)|^{-2s}.
\end{split}
\end{equation*}
Next using \eqref{eq:multcocycle}, we get 
\begin{equation*}
\begin{split}
& E_{k+1/2,s}(\gamma\tau)=\compconj{v_{\vartheta}(\gamma)}^{-2k-1}j(\gamma,\tau)^{k+1/2}|j(\gamma,\tau)|^{2s}M\sum_{g\in\Gamma_0(4)_0\setminus \Gamma_0(4)}\compconj{v_{\vartheta}(g\gamma)}^{2k+1}\Sigma(S^{-1},g\gamma)\\&\hspace{10cm}\times j(S^{-1}g\gamma,\tau)^{-k-1/2}|j(S^{-1}g\gamma,\tau)|^{-2s}\\&=v_{\vartheta}^{2k+1}j(\gamma,\tau)^{k+1/2}|j(\gamma,\tau)|^{2s}M\sum_{g\in\Gamma_0(4)_0\setminus \Gamma_0(4)}\compconj{v_{\vartheta}(g\gamma)}^{2k+1}\Sigma(S^{-1},g\gamma)j(S^{-1}g\gamma,\tau)^{-k-1/2}\\&\hspace{10cm}\times|j(S^{-1}g\gamma,\tau)|^{-2s}\\&=v_{\vartheta}^{2k+1}j(\gamma,\tau)^{k+1/2}|j(\gamma,\tau)|^{2s}E_{k+1/2,s}(\tau),
\end{split}
\end{equation*}
where we used  in third step and absolute convergence in the last step to rearrange the above sum.
\end{proof}
\end{thm}
The modularity transformation of $F_{k+1/2,s}$ follows from Theorem \ref{thm 2.5} and the relation between $E_{k+1/2,s}$ and $F_{k+1/2,s}$. Combining the modularity transformation of $E_{k+1/2,s}$ and $F_{k+1/2,s}$, we obtain the modularity transformation of $H_{k+1/2,s}$. 
\subsection{Review of Zagier's results} 
Let us record a result due to Zagier  which will be used crucially in our calculations.
\theoremstyle{thm}
\begin{prop}\emph{(Zagier, \cite{DZ2})}
For positive integers a and $c,$ let
\begin{equation}
\lambda(a, c)=\left\{\begin{array}{ll}
i^{\frac{1-c}{2}}\left(\frac{a}{c}\right) & \text { if } c \text { is odd, } a \text { even } \\
i^{\frac{a}{2}}\left(\frac{c}{a}\right) & \text { if } a \text { is odd, } c \text { even } \\
0 & \text { otherwise. }
\end{array}\right.
\label{12}
\end{equation}
Define the Gauss sum $\gamma_{c}(n)$ by
\begin{equation}
\gamma_{c}(n):=\frac{1}{\sqrt{c}} \sum_{a=1}^{2 c} \lambda(a, c) e^{-\pi i n \frac{a}{c}}.
\label{13}
\end{equation}
Let $n$ be a non-zero integer and define a Dirichlet series $E_{n}(s)$ by
\begin{equation}
E_{n}(s):=\frac{1}{2} \sum_{c=1 \atop c~\text{odd}}^{\infty} \frac{\gamma_{c}(n)}{c^{s}}+\frac{1}{2} \sum_{c=2 \atop c \text { even}}^{\infty} \frac{\gamma_{c}(n)}{(c / 2)^{s}}
\label{14}
\end{equation}
$($i.e. $E_{n}(s)=\sum a_{m} m^{-s}$ where $a_{m}=\frac{1}{2}\left(\gamma_{m}(n)+\gamma_{2 m}(n)\right)$ when $m$ is odd, and $a_{m}=\frac{1}{2} \gamma_{2 m}(n)$ when
$m \text { is even}).$ Let $K=\mathbb{Q}(\sqrt{n}), D$ be the discriminant of $K, \chi_{d}=\left(\frac{d}{.}\right)$ be the character of $K,$ and $L\left(s, \chi_{d}\right)=\sum \frac{\chi_{D}(n)}{n^{s}}$ be the L-series of $K$ $($if $n$ is a perfect square, then $\chi(m)=1$ for any $m$ and $L(s, \chi)=\zeta(s)) .$ Then if $n \equiv 2,3(\bmod ~4),$ we have
\[
E_{n}(s)=0.
\]
If $n \equiv 0,1(\bmod ~4),$ we have

\begin{equation}
E_{n}(s)=\frac{L\left(s, \chi_{d}\right)}{\zeta(2 s)} \sum_{a, c \geq 1 \atop a c | f} \frac{\mu(a) \chi_{d}(a)}{c^{2 s-1} a^{s}}=\frac{L\left(s, \chi_{d}\right)}{\zeta(2 s)} \frac{T_{s}^{\chi_{d}}(f)}{f^{2 s-1}},
\label{15}
\end{equation}
where $n=df^{2}$ and
\begin{equation}
T_{s}^{\chi}(f)=\sum_{t | f} t^{2 s-1} \sum_{a | t} \frac{\mu(a) \chi(a)}{a^{s}}=\sum_{a | f} \mu(a) \chi(a) a^{s-1} \sigma_{2 s-1}(f/ a).
\label{16}
\end{equation}
Furthermore, we have
\begin{equation}
E_{0}(s)=\frac{\zeta(2 s-1)}{\zeta(2 s)}.
\label{17}
\end{equation}
\label{prop 3.1}
\end{prop}
Let us also define
\begin{equation}
\quad E_{n}^{\text{odd}}(s)=\sum_{c=1 \atop c~\text{odd}}^{\infty} \gamma_{c}(n) c^{-s}
\label{18}
\end{equation}
and
\begin{equation}
E_{n}^{\text {even}}(s)=\sum_{c=1 \atop c~\text {even}}^{\infty} \gamma_{c}(n)(c / 2)^{-s}
\label{19}
\end{equation}
so that
\begin{equation}
E_{n}(s)=\frac{1}{2}\left(E_{n}^{\text{odd}}(s)+E_{n}^{\text {even}}(s)\right)
\label{20}
\end{equation}
\theoremstyle{thm}
\begin{remark}
In \cite{DZ1}, Zagier proved that $E_n(s)$ in Proposition \ref{prop 3.1} can be meromorphically continued to the whole $s-$plane.
\label{remark 2.7}
\end{remark}
\subsection{Meromorphic Continuation of $H_{k+1/2,s}$}
To show that $H_{k+1/2,s}$ has analytic continuation to the whole $s$ plane, we will now calculate the Fourier expansion of $H_{k+1/2,s}$ (which exists because $H_{k+1/2,s}(\tau+1)=H_{k+1/2,s}(\tau)$). The Fourier coefficients may depend on $s$ and Im$(\tau)$. We will then show that the Fourier coefficients are analytic in $s$. The function $H_{k+1/2,s}$ may not be holomorphic in $\tau$ but it satisfies the correct modularity transformation. We begin by calculating the Fourier expansion of $E_{k+1/2,s}$ and $F_{k+1/2,s}$. Observe that 
\begin{equation}
\begin{split}
E_{k+1/2,s}&=\sum_{\substack{m=1\\m~\text{odd}}}^{\infty}\varepsilon_m^{-2k-1}m^{-k-1/2-2s}\sum_{n=-\infty}^{\infty}\left(\frac{n}{m}\right)\left(\tau+\frac{n}{m}\right)^{-k-1/2}\bigg{|}\tau+\frac{n}{m}\bigg{|}^{-2s}
\\&=\sum_{\substack{m=1\\m~\text{odd}}}^{\infty}\varepsilon_m^{-2k-1}m^{-k-1/2-2s}\sum_{n(\text{mod}~m)}\left(\frac{n}{m}\right)\sum_{h\in\mathbb{Z}}\left(\tau+\frac{n}{m}+h\right)^{-k-1/2}\bigg{|}\tau+\frac{n}{m}+h\bigg{|}^{-2s},
\end{split}
\label{21}
\end{equation}
where we replaced $n=mh+k,~0\leq k\leq n-1$ and used the property of Jacobi symbol to replace $\left(\frac{mh+k}{m}\right)$  by $\left(\frac{k}{m}\right)$ in the second step.
By Poisson summation formula,
\begin{equation}
\sum_{h\in\mathbb{Z}}(\tau+h)^{-k-1/2}|\tau+h|^{-2s}=\sum_{h=-\infty}^{\infty}\varrho_h^{k+1/2}(s,v)e^{2\pi ih\tau},~~~~~~~\tau\in\mathbb{H}.
\label{22}
\end{equation}
where 
\begin{equation}
\begin{split}
\varrho_h^{k+1/2}(s,v)&=\int\limits_{iv-\infty}^{iv+\infty}\tau^{-k-1/2}|\tau|^{-2s}e^{-2\pi ih\tau}d\tau\\&=v^{-k+1/2-2s}e^{2\pi hv}\int\limits_{-\infty}^{\infty}(t+i)^{-k-1/2}(t^2+1)^{-s}e^{-2\pi ihtv}dt\\&=v^{-k+1/2-2s}e^{2\pi hv}\int\limits_{-\infty}^{\infty}(t+i)^{-k-1/2-s}(t-i)^{-s}e^{-2\pi ihtv}dt,
\end{split}
\label{23}
\end{equation}
where we got the second expression by substituting $\tau=(t+i)v$.
To simplify $\varrho$, introduce the function
\begin{equation*}
\xi(y ; \alpha, \beta ; t)=\int\limits_{-\infty}^{\infty} e^{-2 \pi i t x}(x+i y)^{-\alpha}(x-i y)^{-\beta} d x
\end{equation*} 
and 
 \begin{equation*}
\Omega(y,\alpha,\beta)\coloneqq \frac{y^{\beta}}{\Gamma(\beta)}\int_{0}^{\infty}e^{-yu}(u+1)^{\alpha-1}u^{\beta-1}du
\end{equation*}
with $y > 0$ and $\alpha,\beta\in\mathbb{C}$ with Re$(\beta) > 0$.
Then we have the following results (cf. Theorem 7.2.5 and Lemma 7.2.4 of \cite{miyake}).
\theoremstyle{thm}
\begin{lemma}
For $(z, \alpha) \in \mathbb{H} \times \mathbb{C},$ we have
\begin{equation*}
\begin{split}
\Omega(z,1-\beta,1-\alpha)=\Omega(z,\alpha,\beta),\\\Omega(z,\alpha,0)=1.
\end{split}
\end{equation*}
\label{lemma 2.8}
\end{lemma} 
\theoremstyle{thm}
\begin{thm}
$(1)$ For each $y>0,$ as a function of $\alpha$ and $\beta,~\xi(y ; \alpha, \beta ; t)$ is analytically continued to a meromorphic function on $\mathbb{C} \times \mathbb{C}$ which is expressed as
\[
\xi(y ; \alpha, \beta ; t)=\left\{\begin{array}{ll}
i^{\beta-\alpha}(2 \pi)^{\alpha} \Gamma(\alpha)^{-1}(2 y)^{-\beta} t^{\alpha-1} e^{-2 \pi y t} \Omega(4 \pi y t ; \alpha, \beta) & (t>0) \\
i^{\beta-\alpha}(2 \pi)^{\alpha+\beta} \Gamma(\alpha)^{-1} \Gamma(\beta)^{-1} \Gamma(\alpha+\beta-1)(4 \pi y)^{1-\alpha-\beta} & (t=0) \\
i^{\beta-\alpha}(2 \pi)^{\beta} \Gamma(\beta)^{-1}(2 y)^{-\alpha}|t|^{\beta-1} e^{-2 \pi y|t|} \Omega(4 \pi y|t| ; \beta, \alpha) & (t<0).
\end{array}\right.
\]
$(2)$ As a function of $\alpha$ and $\beta,~\xi(y ; \alpha, \beta ; t)$ is holomorphic on $\mathbb{C} \times \mathbb{C}$ for $t \neq 0,$ and \\$\Gamma(\alpha+\beta-1)^{-1} \xi(y ; \alpha, \beta ; 0)$ is holomorphic on $\mathbb{C} \times \mathbb{C}.$
\label{thm 2.9}
\end{thm}
Observe that we have
\begin{equation*}
\varrho^{k+1/2}_h(s,v)=v^{-k+1/2-2s}e^{2\pi hv}\xi(1;k+1/2+s,s,hv).
\end{equation*}
Thus Theorem \ref{thm 2.9} implies
\begin{equation}
\varrho_h^{k+1/2}(s,v)=\begin{cases}
(-2i\pi)^{k+1/2}\pi^{s}\Gamma(k+1/2+s)^{-1}h^{k-1/2+s}v^{-s}\Omega(4\pi hv,k+1/2+s,s)&\text{if $h>0$}\\
i^{-k-1/2}(2\pi)\Gamma(k+1/2+s)^{-1}\Gamma(s)^{-1}\Gamma(k-1/2+2s)(2v)^{-k+1/2-2s}&\text{if $h=0$}\\
(2i)^{-k-1/2}\pi^s\Gamma(s)^{-1}(-h)^{s-1}v^{-k-1/2-s}e^{4\pi hv}\Omega(-4\pi hv,s,k+1/2+s)&\text{if $h<0.$}\\
\end{cases}
\label{24}
\end{equation}
Substituting  \eqref{22} in \eqref{21} we get,
\begin{equation}
\begin{split}
E_{k+1/2,s}&=\sum_{\substack{m=1\\m~\text{odd}}}^{\infty}\varepsilon_m^{-2k-1}m^{-k-2s}\sum_{h=-\infty}^{\infty}\Upsilon_m^k(h)\varrho_h^{1+1/2}(s,v)e^{2\pi ih\tau}\\
&=\sum_{h=-\infty}^{\infty}E_{h}^{\text{odd}}(k,k+2s)\varrho_h^{k+1/2}(s,v)e^{2\pi ih\tau},
\end{split}
\label{25}
\end{equation} 
where $\Upsilon_m^k(h)$ is a Gauss sum 
\begin{equation}
\Upsilon_m^k(h)\coloneqq \varepsilon_m^{-2k-1}m^{-1/2}\sum_{n(\text{mod}~m)}\left(\frac{n}{m}\right)e^{2\pi inh/m},~~~~~(m~~\text{odd})
\label{26}
\end{equation}
and $E_{h}^{\text{odd}}(k,k+2s)$ is the Dirichlet series
\begin{equation}
E_{h}^{\text{odd}}(k,s)=\sum_{\substack{m=1\\m~\text{odd}}}^{\infty}\Upsilon_m^k(h)m^{-s}.
\label{27}
\end{equation}
We claim that for $k$ odd, 
\begin{equation}
i^{\frac{1-m}{2}}\left(\frac{2}{m}\right)=\varepsilon_m^{-2k-1}.
\label{28}
\end{equation}
Indeed observe that for $k$ odd, $\varepsilon_m^{-2k-1}=\varepsilon_m$. Then we have that 
\[
i^{\frac{1-m}{2}}\left(\frac{2}{m}\right)  =
  \begin{cases}
                                   1 & \text{if $m~\equiv~1,5(\text{mod}~8)$} \\
                                
i & \text{if $m~\equiv~3,7(\text{mod}~8)$} \\  
  \end{cases}\Bigg{\}}=\varepsilon_m.
\]
Similarly for $k$ even, we have 
\begin{equation}
i^{\frac{1-m}{2}}\left(\frac{-2}{m}\right)=\varepsilon_m^{-2k-1}.
\label{29}
\end{equation}
Now by \eqref{13} we have
\begin{equation*}
\gamma_{m}(h)=\frac{1}{\sqrt{m}} \sum_{n=1}^{2 m} \lambda(n, m) e^{-\pi i h \frac{n}{m}}.
\end{equation*}
Since $\lambda(n,m)=0$ for $n$ odd, thus the sum becomes 
\begin{equation*}
\begin{split}
\gamma_{m}(h)&=\frac{1}{\sqrt{m}} \sum_{n=1}^{m} \lambda(2n, m) e^{-2\pi i h \frac{n}{m}}\\&=\frac{1}{\sqrt{m}}i^{\frac{1-m}{2}}\left(\frac{2}{m}\right) \sum_{n=1}^{m} \left(\frac{n}{m}\right) e^{-2\pi i h \frac{n}{m}}\\&=\Upsilon_m^k(-h),~~~~~(k,m~\text{odd}),
\end{split}
\end{equation*}
where we used \eqref{28}. Similarly $\Upsilon_m^k(h)=\gamma_m(h)$ for $k$ even. Combining these two observations, we have
\begin{equation}
\Upsilon_m^k(h)=\gamma_m((-1)^kh),~~~(m~\text{odd}).
\label{30}
\end{equation} 
In the notation of \ref{18}, we get
\begin{equation}
E_{h}^{\text{odd}}(k,s)=E_{(-1)^kh}^{\text{odd}}(s).
\label{31}
\end{equation}
We will now perform a similar calculation for $F_{k+1/2,s}$. We have
\begin{equation*}
\begin{split}
F_{k+1/2,s}&=\tau^{-k-1/2}|\tau|^{2s}\sum_{\substack{m=1\\m~\text{odd}}}^{\infty}\sum_{\substack{n=-\infty\\\text{gcd}(n,m)=1}}^{\infty}\frac{\left(\frac{n}{m}\right)\varepsilon_m^{-2k-1}}{(-m/4\tau+n)^{k+1/2}|-m/4\tau+n|^{2s}}\\&=\tau^{-k-1/2}|\tau|^{2s}\sum_{\substack{m=1\\m~\text{odd}}}^{\infty}\sum_{\substack{n=-\infty\\\text{gcd}(n,m)=1}}^{\infty}\frac{\left(\frac{n}{m}\right)\varepsilon_m^{-2k-1}}{(-m/4\tau+n)^{k+1/2}|-m/4\tau+n|^{2s}}\\
&=4^{k+1/2+2s}\sum_{\substack{m=1\\m~\text{odd}}}^{\infty}\sum_{\substack{n=-\infty\\\text{gcd}(n,m)=1}}^{\infty}\frac{\text{sign}(n)\left(\frac{n}{m}\right)\varepsilon_m^{-2k-1}}{(-m+4n\tau)^{k+1/2}|-m+4n\tau|^{2s}},
\end{split}
\end{equation*} 
where we used Lemma \ref{lemma:negsqrootcomp} (ii). Replacing $4n$ by $n$ we get
\begin{equation*}
F_{k+1/2,s}=4^{k+1/2+2s}\sum_{\substack{m=1\\m~\text{odd}}}^{\infty}\sum_{\substack{4|n\\\text{gcd}(n,m)=1}}\frac{\text{sign}(n)\left(\frac{n}{m}\right)\varepsilon_m^{-2k-1}}{(n\tau-m)^{k+1/2}|n\tau-m|^{2s}},
\end{equation*} 
where we used the fact that $\left(\frac{4}{m}\right)=1$. Write $m=nh+\ell$. Since $m>0$, we can choose $h>0,~0\leq \ell\leq n-1$ for $n>0$ and $h<0,~n+1\leq \ell\leq 0$ for $n<0$. Since $m=1$ if $n=0$, thus we can write the term for $(m,n)=(1,0)$ separately and then partition the sum as 
\begin{equation*}
\begin{split}
F_{k+1/2,s}&=4^{k+1/2+2s}+\\&4^{k+1/2+2s}\sum_{\substack{n>0\\4|n}}n^{-k-1/2-2s}\sum_{\ell=0}^{n-1}\left(\frac{n}{\ell}\right)\varepsilon_{\ell}^{-2k-1}\sum_{h>0}\left(\tau-\frac{\ell}{n}-h\right)^{-k-1/2}\bigg{|}\tau-\frac{\ell}{n}-h\bigg{|}^{-2s}+\\&4^{k+1/2+2s}\sum_{\substack{n<0\\4|n}}n^{-k-1/2}|n|^{-2s}\sum_{\ell=n+1}^{0}\text{sign}(n)\left(\frac{n}{\ell}\right)\varepsilon_{\ell}^{-2k-1}\sum_{h<0}\left(\tau-\frac{\ell}{n}-h\right)^{-k-1/2}\bigg{|}\tau-\frac{\ell}{n}-h\bigg{|}^{-2s},
\end{split}
\end{equation*}
where we used (c) of Theorem \ref{thm 2.1} $\left(\frac{n}{nh+\ell}\right)=\left(\frac{n}{\ell}\right)$ since $4|n$ and $\varepsilon_{\ell}=\varepsilon_{nh+\ell}$ as $nh+\ell~\equiv~\ell~(\text{mod}~4)$. Now observe that if we replace $n$ by $-n$ and $\ell$ by $-\ell$ in the second sum above, the sum becomes
\begin{equation*}
\begin{split}
& 4^{k+1/2+2s}\sum_{\substack{n>0\\4|n}}(-n)^{-k-1/2}|n|^{-2s}\sum_{\ell=0}^{n-1}\text{sign}(-n)\left(\frac{-n}{-\ell}\right)\varepsilon_{-\ell}^{-2k-1}\sum_{h<0}\left(\tau-\frac{-\ell}{-n}-h\right)^{-k-1/2}\bigg{|}\tau-\frac{-\ell}{-n}-h\bigg{|}^{-2s}\\
&=4^{k+1/2+2s}\sum_{\substack{n>0\\4|n}}n^{-k-1/2-2s}\sum_{\ell=0}^{n-1}(-1)^{-k-1/2}\left(\frac{n}{\ell}\right)\left(\frac{-1}{\ell}\right)i^{-2k-1}\varepsilon_{\ell}^{2k+1}\times\\&\hspace{8.5cm}\sum_{h<0}\left(\tau-\frac{\ell}{n}-h\right)^{-k-1/2}\bigg{|}\tau-\frac{\ell}{n}-h\bigg{|}^{-2s},
\end{split}
\end{equation*}
where we used the fact that $\left(\frac{-n}{-\ell}\right)=\text{sgn}(-n)\left(\frac{n}{\ell}\right)\left(\frac{-1}{\ell}\right)$.
Now we claim that 
\begin{equation}
(-1)^{-k-1/2}i^{-2k-1}\left(\frac{-1}{\ell}\right)\varepsilon_{\ell}^{2k+1}=\varepsilon_{\ell}^{-2k-1}.
\label{32}
\end{equation} 
First suppose $k$ is odd, then $\varepsilon_{\ell}^{2k+1}=\varepsilon_{\ell}^{-1}$. In this case, the left hand side is 
\begin{equation*}
(-i)^{-2k-1}i^{-2k-1}\varepsilon_{\ell}^2\varepsilon_{\ell}^{-1}=\varepsilon_{\ell}=\varepsilon_{\ell}^{-2k-1}.
\end{equation*} 
Here we used the fact that $(-1)^{-k-1/2}=(-i)^{-2k-1}$ in the principal branch of square root. Similarly if $k$ is even, then $\varepsilon_{\ell}^{2k+1}=\varepsilon_{\ell}$. In this case, the left hand side is
\begin{equation*}
(-i)^{-2k-1}i^{-2k-1}\varepsilon_{\ell}^2\varepsilon_{\ell}=\varepsilon_{\ell}^{-1}=\varepsilon_{\ell}^{-2k-1}.
\end{equation*}
Using \eqref{32} and the calculations above, we get
\begin{equation*}
\begin{split}
F_{k+1/2,s}&=4^{k+1/2+2s}+\\&4^{k+1/2+2s}\sum_{\substack{n>0\\4|n}}n^{-k-1/2-2s}\sum_{\ell=0}^{n-1}\left(\frac{n}{\ell}\right)\varepsilon_{\ell}^{-2k-1}\sum_{h=-\infty}^{\infty}\left(\tau-\frac{\ell}{n}-h\right)^{-k-1/2}\bigg{|}\tau-\frac{\ell}{n}-h\bigg{|}^{-2s}.
\end{split}
\end{equation*}
Changing $h$ to $-h$ and using the Poisson summation formula of \eqref{22}, we get
\begin{equation*}
\begin{split}
F_{k+1/2,s}&=4^{k+1/2+2s}+\\&4^{k+1/2+2s}\sum_{\substack{n>0\\4|n}}n^{-k-1/2-2s}\sum_{\ell=0}^{n-1}\left(\frac{n}{\ell}\right)\varepsilon_{\ell}^{-2k-1}\sum_{h=-\infty}^{\infty}\left(\tau-\frac{\ell}{n}+h\right)^{-k-1/2}\bigg{|}\tau-\frac{\ell}{n}+h\bigg{|}^{-2s}\\&=4^{k+1/2+2s}+\\&~~~~~~~~~~~~4^{k+1/2+2s}\sum_{\substack{n>0\\4|n}}n^{-k-1/2-2s}\sum_{\ell=0}^{n-1}\left(\frac{n}{\ell}\right)\varepsilon_{\ell}^{-2k-1}\sum_{h=-\infty}^{\infty}\varrho_h^{k+1/2}(s,v)e^{-2\pi i\ell h/n}e^{2\pi ih\tau}\\&=4^{k+1/2+2s}+\\&~~~~\sum_{h=-\infty}^{\infty}\left[4^{k+1/2+2s}\sum_{\substack{n>0\\4|n}}\left(\sum_{\ell=0}^{n-1}n^{-1/2}\left(\frac{n}{\ell}\right)\varepsilon_{\ell}^{-2k-1}e^{-2\pi i\ell h/n}\right)n^{-k-2s}\right]\varrho_h^{k+1/2}(s,v)e^{2\pi ih\tau}\\&=4^{k+1/2+2s}(-i)^{k}+\\&~~~~\sum_{h=-\infty}^{\infty}\left[\sqrt{2}\sum_{\substack{n=1\\n~\text{even}}}^{\infty}\left(\sum_{\ell=0}^{2n-1}n^{-1/2}\left(\frac{2n}{\ell}\right)\varepsilon_{\ell}^{-2k-1}e^{-2\pi i\ell h/2n}\right)(n/2)^{-k-2s}\right]\varrho_h^{k+1/2}(s,v)e^{2\pi ih\tau},
\end{split}
\end{equation*}
where we replaced $-(-i)^{-2k-1}=(-i)^{k}$. Next, observe that 
\begin{equation}
\begin{split}
\sqrt{2}\left(\frac{2}{\ell}\right)i\varepsilon_{\ell}^{-1}=(1+i)i^{\ell/2}\\\sqrt{2}\left(\frac{2}{\ell}\right)\varepsilon_{\ell}^{-1}=(1-i)i^{\ell/2}.
\end{split}
\label{33}
\end{equation}
This is verified easily for $\ell~\equiv~1,3,5,7(\text{mod}~8)$. 
We then note that for $k$ odd, replacing $\varepsilon_{\ell}^{-2k-1}$ by $\varepsilon_{\ell}$ and changing $\ell$ to $-\ell$ in the inner most sum in the above expression we get
\begin{equation*}
\begin{split}
\sqrt{2}n^{-1/2}\sum_{\ell=0}^{2n-1}\left(\frac{2n}{\ell}\right)\varepsilon_{\ell}^{-2k-1}e^{-2\pi i\ell h/2n}&=\sqrt{2}n^{-1/2}\sum_{\ell=0}^{2n-1}\left(\frac{2n}{\ell}\right)\sqrt{2}i\varepsilon_{\ell}^{-1}e^{2\pi i\ell h/2n}\\&=(1+i)n^{-1/2}\sum_{\ell=0}^{2n-1}\left(\frac{n}{\ell}\right)i^{\ell/2}e^{2\pi i\ell h/2n}\\&=(1+i)\gamma_n(-h),
\end{split}
\end{equation*}
where $\gamma_n(h)$ is defined in Proposition \ref{prop 3.1}. Here we used the fact that $\varepsilon_{-\ell}=i\varepsilon_{\ell}^{-1}$ and $\left(\frac{n}{-1}\right)=\text{sign}(n)=1$ and \ref{33}.  Similarly for $k$, replacing $\varepsilon_{\ell}^{-2k-1}$ by $\varepsilon_{\ell}^{-1}$ the inner most sum becomes 
\begin{equation*}
\begin{split}
\sqrt{2}n^{-1/2}\sum_{\ell=0}^{2n-1}\left(\frac{2n}{\ell}\right)\varepsilon_{\ell}^{-1}e^{-2\pi i\ell h/2n}&=(1-i)n^{-1/2}\sum_{\ell=0}^{2n-1}\left(\frac{n}{\ell}\right)i^{\ell/2}e^{-2\pi i\ell h/2n}\\&=(1-i)\gamma_n(h).
\end{split}
\end{equation*} 
Combining the above two observations, we get 
\begin{equation}
\sqrt{2}n^{-1/2}\sum_{\ell=0}^{2n-1}\left(\frac{2n}{\ell}\right)\varepsilon_{\ell}^{-1}e^{-2\pi i\ell h/2n}=(1+i^{2k-1})\gamma_n((-1)^kh).
\label{34}
\end{equation}
Using \ref{34}, we get 
\begin{equation}
F_{k+1/2,s}=2^{4s+2k+1}+(1+i^{2k-1})\sum_{h=-\infty}^{\infty}E_{h}^{\text{even}}(k,k+2s)\varrho_h^{k+1/2}(s,v)e^{2\pi ih\tau},
\label{35}
\end{equation} 
where
\begin{equation}
E_{h}^{\text{even}}(k,s)=\sum_{\substack{m=1\\m~\text{even}}}^{\infty}\Upsilon_m^k(h)(m/2)^{-s},
\label{36}
\end{equation}
and the Gauss sum $\Upsilon_m^k(h)=\gamma_m((-1)^kh).$ In the notation of Proposition \ref{prop 3.1}
\begin{equation}
E^{\text{even}}_h(k,s)=E^{\text{even}}_{(-1)^kh}(s).
\label{37}
\end{equation}
We now have the following proposition.
\begin{prop}
Let $k\in\mathbb{Z}$ and $s\in\mathbb{C}$ with Re$(s)>\frac{3}{4}-\frac{k}{2}$, we have 
\begin{equation}
H_{k+1/2,s}(\tau)\coloneqq \frac{\zeta(1-2k)}{2^{-4s}}i^{2k+1}+\frac{(1+i^{2k+1})\zeta(1-2k)}{2^{2k}}\sum_{h=-\infty}^{\infty}E_{h}(k,k+2s)\varrho_h^{k+1/2}(s,v)e^{2\pi ih\tau},
\label{38}
\end{equation}
where in the notation of \eqref{20} 
\begin{equation}
E_{h}(k,s)\coloneqq \frac{1}{2}\left(E_{h}^{\text{even}}(k,s)+E_{h}^{\text{odd}}(k,s)\right)=E_{(-1)^kh}(s).
\label{39}
\end{equation}
In particular $H_{k+1/2,s}$ can be meromorphically continued to the whole $s-$plane.
\label{prop 3.10}
\begin{proof}
Combining \eqref{25} and \eqref{35} and using \eqref{31} and \eqref{37} immediately gives the Fourier series expansion of $H_{k+1/2,s}$. Next, by Remark \ref{remark 2.7} and Theorem \ref{thm 2.9}, $E_{(-1)^kh}(k+2s)\varrho_h^{k+1/2}(s)$ can be meromorphically continued to the whole $s-$plane. Thus $H_{k+1/2,s}$ can be meromorphically continued to the whole $s-$plane using the above Fourier expansion. 
\end{proof} 
\end{prop}
Thus we can define $H_{3/2}=\lim\limits_{s\to 0}H_{3/2,0}$.
In subsequent computations, we will need the functional equation for Dirichlet $L-$function \cite{Ram murty}, so we record it here for completeness.
\begin{equation}
L(s,\chi)=\omega_{\chi}.\begin{cases}
                                   N^{-s+1/2}\pi^{s-1/2}\Gamma(1-s/2)\Gamma\left(\frac{s+1}{2}\right)^{-1}L(1-s,\chi) & \text{if $\chi(-1)=-1$,} \\
                                   N^{-s+1/2}\pi^{s-1/2}\Gamma(s/2)^{-1}\Gamma\left(\frac{1-s}{2}\right)L(1-s,\chi) & \text{if $\chi(-1)=1$,} \\
  \end{cases}
\label{40}  
\end{equation}
where $\chi$ is a Dirichlet character modulo $N$ and 
\begin{equation*}
\omega_{\chi}=\begin{cases}\frac{\tau_{\chi}}{i\sqrt{N}} & \text{if $\chi(-1)=-1$,} \\
                                   \frac{\tau_{\chi}}{\sqrt{N}} & \text{if $\chi(-1)=1$,} \\
  \end{cases}
\end{equation*}
where $\tau_{\chi}=\sum\limits_{b=0}^{N-1}\chi(b)e^{2\pi ib/N}$ is the Gauss sum associated to the character $\chi$.
\section{Proof of Theorem \ref{thm 1.3}}
\label{sec 4}
By Proposition \ref{prop 3.10}, we have
\begin{equation}
\begin{split}
H_{3/2,s}(\tau)&= \frac{\zeta(-1)}{2^{-4s}}+\frac{(1-i)\zeta(-1)}{4}\sum_{h=-\infty}^{\infty}E_{-h}(1+2s)\varrho_h^{3/2}(s,v)e^{2\pi ih\tau}\\
&=-\frac{4^{4s}}{12}-\frac{(1-i)}{48}\sum_{h=-\infty}^{\infty}E_{-h}(1+2s)\varrho_h^{3/2}(s,v)e^{2\pi ih\tau},
\end{split}
\label{41}
\end{equation}
where we used $\zeta(-1)=-1/2.$ If we write 
\begin{equation}
H_{3/2,s}=\sum_{h=-\infty}^{\infty}\alpha_h(s,v)q^h,
\label{42}
\end{equation}
then we would like to show that 
\begin{equation}
\alpha_h(0,v)=\begin{cases}
-\frac{1}{12}+\frac{1}{8\pi\sqrt{v}}&\text{if $h=0$}\\
H(h)&\text{if $h>0$}\\
\frac{f}{4\sqrt{\pi}}\Gamma\left(-\frac{1}{2},4\pi f^2v\right)&\text{if $h=-f^2,~f>0$}\\
0&\text{otherwise.}\\
\end{cases}
\label{43}
\end{equation}
By Proposition \ref{prop 3.1} we have 
\begin{equation}
 E_{h}(s)  =
  \begin{cases}
                                   0 & \text{if $h~\equiv~2,3(\text{mod}~4)$} \\
                                   \frac{L(s,\chi_d)}{\zeta(2s)}.\frac{T^{\chi_d}_s(f)}{f^{2s-1}} & \text{if $h~\equiv~0,1(\text{mod}~4), h\neq 0.$} \\
\frac{\zeta(2s-1)}{\zeta(2s)} & \text{if $h=0,$} \\  
  \end{cases}
\label{44}  
\end{equation}
where $\chi_d$ and $T^{\chi_d}_s(f)$ is defined in Eqs. \eqref{6} and \eqref{7}.
Using the integral representation of zeta function \cite[Theorem 5.2.2]{Ram murty}, we have
\begin{equation}
\lim_{s\to 0}\frac{\zeta(4s+1)}{\Gamma(s)}=\lim_{s\to 0}\frac{4s\zeta(4s+1)\Gamma\left(\frac{4s+1}{2}\right)}{4s\Gamma(s)\Gamma\left(\frac{4s+1}{2}\right)}=\frac{1}{4}\frac{\sqrt{\pi}}{\sqrt{\pi}}=\frac{1}{4},
\label{45}
\end{equation}
where we used $z\Gamma(z)=\Gamma(z+1)$.
Then using  \eqref{24} we have
\begin{equation*}
\begin{split}
\lim_{s\to 0}E_0(1+2s)\varrho_0^{3/2}(0,v)&=i^{-3/2}(2\pi)\Gamma(3/2)^{-1}\frac{1}{4}\Gamma(1/2)\zeta(2)^{-1}v^{-1/2}2^{-1/2}
\\&=\frac{6i^{-3/2}}{\pi\sqrt{2}\sqrt{v}}.
\end{split}
\end{equation*} 
Now using \eqref{41}, we get
\begin{equation*}
\alpha_0(0,v)=-\frac{1}{12}-\frac{(1-i)}{48}\frac{6i^{-3/2}}{\pi\sqrt{2}\sqrt{v}}=-\frac{1}{12}+\frac{1}{8\pi\sqrt{v}},
\end{equation*}
where we used \eqref{33} for $\ell=-3$. Now for $h~\equiv~0,3(\text{mod}~4),~h>0$ we have
\begin{equation*}
E_{-h}(k)=\frac{L(k,\chi_d)}{\zeta(2k)}.\frac{T^{\chi_d}_k(f)}{f^{2k-1}}.
\end{equation*}
If we write $-h=df^2$ then $d<0$. This implies that $\chi_d$ is an odd character mod $-d$. Using the functional equation for $L(s,\chi)$ from \eqref{40}, we get
\begin{equation*}
\begin{split}
L(k,\chi_d)=\omega_{\chi_d}(-d)^{-k+1/2}\pi^{k-1/2}\Gamma(1-k/2)\Gamma\left(\frac{k+1}{2}\right)^{-1}L(1-k,\chi_d).
\end{split}
\end{equation*}
One can easily evaluate $\tau_{\chi_d}$ using the quadratic reciprocity law for Kronecker symbol. It turns out that (cf. Theorem 9.17 of \cite{monty}) 
\begin{equation*}
\tau_{\chi_d}=\begin{cases}
\sqrt{d}&\text{if $d>0$}\\
i\sqrt{-d}&\text{if $d<0.$}
\end{cases}
\end{equation*}
This means that $\omega_{\chi_d}=1.$
Thus using  \eqref{40} we get for $k$ odd
\begin{equation*}
\begin{split}
E_{-h}(k)&=\frac{1}{\zeta(2k)}(-df^2)^{-k+1/2}\pi^{k-1/2}\Gamma(1-k/2)\Gamma\left(\frac{k+1}{2}\right)^{-1}L(1-k,\chi_d)T^{\chi_d}_k(f)\\&=\frac{1}{\zeta(2k)}h^{-k+1/2}\pi^{k-1/2}\Gamma(1-k/2)\Gamma\left(\frac{k+1}{2}\right)^{-1}H(k,h).
\end{split}
\end{equation*} 
In particular 
\begin{equation*}
E_{-h}(1)=\frac{6}{\pi^2}h^{-1/2}\pi H(h).
\end{equation*} 
Using \eqref{24} for $h>0$ we get
\begin{equation*}
\begin{split}
\varrho_h^{3/2}(0,v)&=i^{-3/2}\pi^{3/2}2^{3/2}\Gamma(3/2)^{-1}h^{1/2}\Omega(4\pi hv,3/2,0)\\&=i^{-3/2}\frac{8\pi h^{1/2}}{\sqrt{2}},
\end{split}
\end{equation*}
where we used Lemma \ref{lemma 2.8}. We have for $h>0$
\begin{equation*}
\alpha_h(0,v)=-\frac{1-i}{48}E_{-h}(1)\varrho_h^{3/2}(0,v)=H(h),
\end{equation*}
where again we used \eqref{33}. Next suppose $h<0$ but $-h$ is not a square, then $-h=df^2$ implies $d\neq 1$. Clearly $\chi_d$ is non-principal character. Thus $L(s,\chi_d)$ is holomorphic at $s=1$. This in turn means that $E_{-h}(s)$ is holomorphic  at $s=1$ (cf. Exercise 2.3.4 of \cite{Ram murty}). Since $\frac{1}{\Gamma(0)}=0$ we have $E_{-h}(1)\varrho_h^{3/2}(0,v)=0$, thus $\alpha_h(0,v)=0$ for $h<0$ and $h$ not a square. Finally, if $h=-f^2,~f>0$, then 
\[
E_{-h}(s)=\frac{\zeta(s)}{\zeta(2s)}\frac{T^{\chi_1}_s(f)}{f^{2s-1}},
\]
since $L(s,\chi_1)=\zeta(s)$. Using M\"{o}bius inversion formula, it is easy to see that 
\[
\frac{T^{\chi_1}_1(f)}{f}=1.
\]
Next observe that for $h<0$ we have
\begin{equation*}
v^{\alpha-2}\Omega(-4\pi hv,\alpha-1,1)=-4\pi hv^{\alpha-1}\int_{0}^{\infty}e^{4\pi hvx}(x+1)^{\alpha-2}dx.
\end{equation*}
Substituting $t=x+1$ we get 
\begin{equation*}
v^{\alpha-2}\Omega(-4\pi hv,\alpha-1,1)=-4\pi he^{-4\pi hv}v\int_{1}^{\infty}e^{4\pi hvt}(vt)^{\alpha-2}dt.
\end{equation*}
Again substituting $w=-4\pi hvt$ we get
\begin{equation*}
v^{\alpha-2}\Omega(-4\pi hv,\alpha-1,1)=(-4\pi h)^{2-\alpha}e^{-4\pi hv}\Gamma(\alpha-1,-4\pi hv).
\end{equation*}
Thus we have
\begin{equation}
\begin{split}
v^{-3/2}e^{4\pi hv}\Omega(-4\pi hv,0,3/2)=v^{-3/2}e^{4\pi hv}\Omega(-4\pi hv,-1/2,1)=(-4\pi h)^{3/2}\Gamma\left(-\frac{1}{2},-4\pi hv\right),
\end{split}
\label{46}
\end{equation}
where we used Lemma \ref{lemma 2.8}. This implies that for $h=-f^2,~f>0$ we have
\begin{equation*}
\begin{split}
E_{-h}(1)\varrho^{3/2}_h(0,v)&=2i^{-3/2}2^{-3/2}\zeta(2)^{-1}\frac{1}{4}\left(\frac{1}{-h}\right)v^{-3/2}e^{4\pi hv}\Omega(-4\pi hv,0,3/2)\\&=\frac{12i^{-3/2}f}{\sqrt{2}\sqrt{\pi}}\Gamma\left(-\frac{1}{2},4\pi f^2v\right),
\end{split}
\end{equation*}
where we used \eqref{45} and \eqref{46}. Finally we get
\[
\alpha_h(0,v)=-\frac{(1-i)}{48}\frac{12i^{-3/2}f}{\sqrt{2}\sqrt{\pi}}\Gamma\left(-\frac{1}{2},4\pi f^2v\right)=\frac{f}{4\sqrt{\pi}}\Gamma\left(-\frac{1}{2},4\pi f^2v\right).
\]
Lastly, computing the shadow of $H_{3/2}(\tau)=H_{3/2,0}(\tau)$ is a trivial exercise using \eqref{4}. Indeed with $2-k=3/2$ we have
\[
\xi_{3/2}(H_{3/2})=-\frac{1}{16\pi}-\frac{(4\pi)^{-1/2}}{4\sqrt{\pi}}\sum_{n=1}^{\infty}(n^2)^{-1/2}nq^{n^2}=-\frac{1}{16}\left[1+2\sum_{n=1}^{\infty}q^{n^2}\right]=-\frac{1}{16}\Theta(\tau).
\]
Since $\xi$ annihilates holomorphic functions, we have that $\Delta_{3/2}(H_{3/2})=-\xi_{1/2}(\xi_{3/2}(H_{3/2}))=0.$ The growth condition follows from the Fourier expansion. This completes the proof of Zagier's theorem.
\bibliographystyle{plain}
\bibliography{mybib}

\end{document}